\documentclass[10pt]{article}
\usepackage{mathrsfs}
\usepackage{stmaryrd}
\usepackage{amsfonts,amsmath,amssymb,amscd}
\usepackage[colorlinks,linkcolor=red,anchorcolor=blue,citecolor=green]{hyperref}
\usepackage{shadow}
\usepackage{graphics,graphicx}
\usepackage{color}
\usepackage{enumerate}
\usepackage{longtable}
\usepackage{appendix}
\allowdisplaybreaks

\newtheorem{thm}{Theorem}[section]
\newtheorem{lem}[thm]{Lemma}

\newtheorem{remark}[thm]{Remark}

\newtheorem{conj}[thm]{Conjecture}

\newtheorem{ass}[thm]{Assumption}
\def\qed{\hfill \rule{4pt}{7pt}}

\def\pf{\noindent {\it{Proof.}\hskip 2pt}}

\begin{document}
\begin{center}
{{\large\bf Connective Constants on Cayley Graphs}}
\end{center}

\begin{center}
Song He, Xiang Kai-Nan and Zhu Song-Chao-Hao
\vskip 1mm
\footnotesize{School of Mathematical Sciences, LPMC, Nankai University}\\
\footnotesize{Tianjin City, 300071, P. R. China}\\
\footnotesize{Emails: songhe@mail.nankai.edu.cn} (Song)\\
\footnotesize{~~~~~~~~~~~~kainanxiang@nankai.edu.cn} (Xiang)\\
\footnotesize{~~~~~~~~~~~~zsch61@qq.com} (Zhu)
\footnote{The project is supported partially by CNNSF (No. 11271204).

   {\it \ \ MSC2010 subject classifications}. 05C30, 82B20, 60K35.

   {\it \ \ Key words and phrases}. Connective constant, locality, Cayley graph.}
\end{center}

\begin{abstract}
For a transitive infinite connected graph $G$, let $\mu(G)$ be its connective constant.
Denote by $\mathbf{\cal G}$ the set of Cayley graphs for finitely generated infinite groups with an infinite-order generator which is independent
of other generators. Assume $G\in\mathbf{\cal G}$ is a Cayley graph of a finitely presented group, and
Cayley graph sequence $\{G_n\}_{n=1}^{\infty}\subset \mathbf{\cal G}$ converges locally to $G.$
Then $\mu(G_n)$ converges to $\mu(G)$ as $n\rightarrow\infty.$ This confirms partially a conjecture raised by Benjamini [2013. {\it Coarse geometry and randomness.} Lect. Notes Math. {\bf 2100}. Springer.] that connective constant is continuous with respect to local convergence of infinite transitive connected graphs.
\end{abstract}

\section{Introduction}
For a locally finite, connected infinite graph $G=(V, E),$ a {\it self-avoiding walk} (SAW) on it is a path that visits each vertex at most one time.
SAW was first introduced by Flory \cite{Flo} in the setting of long-chain polymers in chemistry, and its critical behavior has been received much attention
by mathematicians and physicists (\cite{Bau-Dum-Goo}, \cite{Mad-Sla}).

Let $c_n(v)$ be the number of $n$-step SAWs on $G$ with an initial vertex $v.$
Define $\mu (G)=\lim\limits_{n\rightarrow\infty}c_n(v)^{\frac{1}{n}}$ if it exists and does not depend on $v.$ Call $\mu(G)$ the {\it connective constant} of $G.$ Recall from Hammersley \cite{Ham}, $\mu (G)\in [1,\infty)$ is well-defined for quasi-transitive $G.$
When $G$ is transitive, $c_n(v)$ is independent of $v$ and denote it by $c_n.$
Connective constants are exactly known only for few graphs. For example, $\mu(\mathbb{Z}^2)$ is unknown. And for the hexagonal lattice $\mathbb{H}$
in a plane, $\mu(\mathbb{H})=\sqrt{2+\sqrt{2}}$ was proven by Duminil-Copin and Smirnov \cite{Dum-Co-Smi} by exploiting the construction of an observable with some properties on discrete holomorphicity and the bridge decomposition introduced in Hammersley and Welsh \cite{Ham-Wel}. This is a very significant recent result.

To continue, assume $G$ is transitive. For a sequence $\{G_n\}_{n=1}^{\infty}$ of transitive graphs, say it converges locally to $G,$ if
for any natural number $r,$ $B_{G_n}(x_n,r)$ is isomorphic to $B_G(x,r)$ when $n$ is large enough. Here for a graph $H$ and its vertex $v,$
$B_H(v,r)$ is the ball in $H$ with radius $r$ and center $v;$ and $x$ (resp. $x_n$) is an arbitrary vertex of $G$ (resp. $G_n$). Recall from Benjamini \cite{Ben} Chapter 4 the following

\begin{conj}\label{conjecture1}
Connective constant $\mu(G)$ is continuous with respect to local convergence of infinite transitive connected graphs $G.$
\end{conj}

Conjecture \ref{conjecture1} is the SAW case for the locality conjecture of
critical parameters in physical systems. And for percolation, the parameter in question is critical probability; while for Ising model, it is critical temperature.
It is important to understand whether critical parameter is locally or globally determined by the geometry of graphs. For the related locality conjecture,
see \cite{Ben-Nac-Per}, \cite{Mar-Tas} and references therein in the percolation case, and \cite{Bod}, \cite{LiZY2} and \cite{Cim-Dum} in the Ising model
setting.

Recall connective constant was studied extensively by Grimmett and Li \cite{Gri-Li1}-\cite{Gri-Li3} recently. And Li \cite{LiZY1} proved Conjecture \ref{conjecture1} for Cayley graphs under some conditions. In the following we describe briefly the result of \cite{LiZY1}.

To begin, let $G=(\Gamma,S)$ be an infinite Cayley graph of a finitely generated group $\Gamma$ with the following finite generating set
$$S=\left\{t_1,\cdots, t_p\right\};$$
where edge set of $G$ is $\left\{(g,gs);\ g\in\Gamma,\ s\ \mbox{or}\ s^{-1}\ \mbox{is in}\ S\right\}.$ Suppose $\Gamma$ has a presentation $\Gamma=\langle S\vert R\rangle$ with $R$ being the relator set. Let $G_m$ be the Cayley graph obtained from $G$ by adding more relators. Here relator means a word of generators that is identified with the identity element of the group, namely a cycle of the Cayley graph. In other words,
$$
G_m=(\Gamma_m,S),\ \Gamma_m=\langle S\vert R_m\rangle, R\subseteq R_m. \eqno{(1.1)}
$$
Define the relative girth $\widetilde{g}_m$ of $G_m$ with respect to $G$ as the minimum length of cycles in $G_m$ but not in $G$
if such circles exist, and otherwise let $\widetilde{g}_m=\infty.$

\begin{ass}\label{Assumption1}
${\bf (i)}$ $\Gamma$ is a finitely generated infinite group, namely $S$ is finite. ${\bf (ii)}$ Every $G_m$ is infinite and connected.
${\bf (iii)}$ $\lim\limits_{m\rightarrow\infty}\widetilde{g}_m=\infty.$
\end{ass}

Note that $G_m$ converges locally to $G$ under  Assumption \ref{Assumption1}.
Associate a linear equation system to each Cayley graph $G_m$ with unknown variables $(\alpha_1,\cdots,\alpha_p).$ Due to each relator in $R_m$ of $\Gamma_m$ is a word consisting of elements in $S,$ we view every relator as a monomial with unknown variables $t_1,\cdots,t_p,$ and construct the following linear equation for any relator in $\Gamma_m:$  When degrees of $t_1,\cdots,t_p$ in a relator are $u_1,\cdots,u_p,$ the constructed equation is $\sum\limits_{j=1}^p\alpha_ju_j=0.$ Denote by
$C_m=\left(u^{(i)}_j\right)_{i,j}$ the corresponding coefficient matrix, where index $(i)$ is used to distinguish different relators. Let $C\alpha=0$ be the linear equation system consisting of all linear equations $C_m\alpha=0,\forall m\geq 1,$ where $\alpha$ is the column vector with the $j$th component being $\alpha_j.$

Write $r(C)$ for the rank of the matrix $C.$
Then Li \cite{LiZY1} Theorem 2 reads as follows: For $\{G_m\}_{m=1}^\infty$ specified in (1.1), $\displaystyle \lim_{m\rightarrow\infty}\mu(G_m)=\mu(G)$
under Assumption \ref{Assumption1} and
$$p>r(C).\eqno{(1.2)}$$

Let $\mathbf{\cal G}$ be the set of Cayley graphs for finitely generated infinite groups with an infinite-order generator which is independent of other generators.
In this paper, we prove the following

\begin{thm}\label{Thm1}
Let $G\in\mathbf{\cal G}$ be a Cayley graph corresponding to a finitely presented group $\Gamma.$
Then for any Cayley graph sequence $\{G_n\}_{n=1}^{\infty}\subset \mathbf{\cal G}$ converging locally to $G,$
$\lim\limits_{n\rightarrow\infty}\mu(G_n)=\mu(G).$
\end{thm}
\vskip 2mm

\begin{remark}\label{rem1}
For a finitely generated infinite group, there may not be an element of infinite-order in general. The assumption that there is an
infinite-order generator independent of other generators ensures the existence of a nontrivial ``invariant" antisymmetric edge-weights on Cayley graphs
and the validity of (2.2). See proof of Lemma \ref{lem1}.

Comparing with \cite{LiZY1}, difference in
proving Lemma \ref{lem1} is that it is unnecessary to define an ``invariant" antisymmetric edge-weight function such that the edge-weight function is nontrivial
restricted to cycles, and the edge-weight sum along any (directed) cycle is zero. Thus (1.2) is unnecessary for Theorem \ref{Thm1} to hold.

Why do we assume $\Gamma$ is finitely presented? It lies in that we need $\Gamma$ is a quotient group of a free group by a finitely
generated normal subgroup, and this is a key point to prove Lemma \ref{lem2}. Hence we do not assume each $G_n$ is a quotient of $G$ in Theorem
\ref{Thm1}, which differs from that of \cite{LiZY1}.

Note Lemmas \ref{lem1} and \ref{lem2} play important roles in proof of Theorem \ref{Thm1}. It is challenging to remove technical condition that $\Gamma$
is finitely presented, and $G$ and every $G_n$ are in $\mathbf{\cal G}.$
\end{remark}

\section{Proof of Theorem \ref{Thm1}}
We firstly prove Lemma \ref{lem1} on some kind of localities for connective constants based on \cite{LiZY1} and some new insights. Then we verify Lemma \ref{lem2} on marked groups, which is an interesting extension of the related version of marked abelian groups in \cite{Mar-Tas}. Finally, by Lemmas \ref{lem1} and \ref{lem2}, we can prove Theorem \ref{Thm1} by reduction to absurdity.

Let $\{H_m\}_{m=1}^{\infty}\subset\mathbf{\cal G}$ be a sequence of Cayley graphs with generating set sequence $\{S_m\}_{m=1}^\infty,$
and $H\in\mathbf{\cal G}$ a Cayley graph with a generating set $S.$
Assume every $H_m$ is a quotient graph of $H,$ and relative girth $\widehat{g}_m$ of $H_m$ with respect to $H$ tends to infinity
as $m\rightarrow\infty.$

\begin{lem}\label{lem1}
For $H_m$ and $H$ specified above, $\lim\limits_{m\rightarrow\infty}\mu(H_m)=\mu(H).$
\end{lem}
\pf
\noindent {\bf Step 1.} {\it Definitions: bridge and half-space walk.}

Let $S_m=\left\{s_1,\cdots,s_{\ell_m}\right\}$ be the finite generating set for Cayley graph $H_m,$
and $s_1$ be of infinite-order and independent of other generators $s_j$ with $2\leq j\leq \ell_m.$
For any directed edge $(x,y)$ of $H_m,$ endow it with a weight as follows:
\[
w(x,y)=\left\{\begin{array}{cl}
      0 &{\rm if}\ x^{-1}y\ \mbox{or}\ y^{-1}x\ \mbox{is in}\ S_m\setminus\left\{s_1\right\},\\
      1 &{\rm if}\ x^{-1}y=s_1\in S_m,\\
      -1 &{\rm if}\ x^{-1}y=s_1^{-1}\in S_m.
   \end{array}
\right.
\]
For any $t\in H_m,$ let $\phi_t$ be the following automorphism of $H_m:$
$x\in H_m\rightarrow tx\in H_m.$
Clearly,
 $$w(x,y)=-w(y,x)\ \mbox{and}\ w(\phi_t(x,y))=w(x,y)\ \mbox{for any directed edge}\ (x,y)\ \mbox{of}\ H_m.\eqno{(2.1)}$$

For any $n\geq 1$ and $n$-step SAW $\omega=(\omega(s))_{0\leq s\leq n}$ of $H_m$ starting at a vertex $a\in H_m,$ the height $h_s(\omega)$ of $\omega(s)$ in $\omega$ is $0$ when $s=0$ and $\sum\limits_{i=1}^sw(\omega(i-1),\omega(i))$ when $s\geq 1.$  Call $\omega$ a {\it bridge} if
$$h_0(\omega)<h_s(\omega)\leq h_n(\omega),\ 1\leq s\leq n;$$
and a {\it half-space walk} if $h_0(\omega)<h_s(\omega),\ 1\leq s\leq n.$ The {\it span} of $\omega$ is $\max\limits_{0\leq s\leq n}h_s(\omega)-\min\limits_{0\leq s\leq n}h_s(\omega).$

Denote the number of $n$-step half-space walks (resp. bridges) starting at $a$ and having span $A$  by $h_{n,A}(a)$ (resp. $b_{n,A}(a)$). By (2.1),
both $h_{n,A}(a)$ and $b_{n,A}(a)$ do not depend on $a.$ Hence write $h_{n,A}$ and $b_{n,A}$ for $h_{n,A}(a)$ and $b_{n,A}(a)$ respectively.
And the number $h_n$ (resp. $b_n$) of $n$-step half-space walks (resp. bridges) starting at any fixed vertex is
$$h_n=\sum\limits_{A=1}^nh_{n,A}\ \left(\mbox{resp.}\ b_n=\sum\limits_{A=1}^nb_{n,A}\right).$$
{\it Convention.} A single point is called a $0$-step half-space walk and a $0$-step bridge. And $h_0=b_0=1.$\\

\noindent{\bf Step 2.} {\it For any $N\geq 1,$ $h_N\leq \sum\limits_{A=1}^{N}P_D(A)b_{N,A}\leq P_D(N)b_N.$ Here $P_D(A)$ is the number of ways to write  $A=A_1+\cdots +A_k$ with $A_1>\cdots>A_k$ being natural numbers.}

{\it Indeed}, let $\omega=(\omega(s))_{0\leq s\leq N}$ be an $N$-step SAW starting from $a\in H_m$ and $n_0=0;$ and define recursively $A_j(\omega)$ and $n_j(\omega)$ for $j=1,2,\cdots$ as follows:
 $$A_j=\max\limits_{n_{j-1}\leq s\leq N}(-1)^j(h_{n_{j-1}}(\omega)-h_s(\omega)),\ n_j=\max\left\{n_{j-1}\leq s\leq N\ \left\vert\
   (-1)^j(h_{n_{j-1}}(\omega)-h_s(\omega))=A_j\right.\right\};$$
and this recursion is terminated at the smallest $k$ with $n_k=N.$ Then $A_{j}$ is the span of SAW $(\omega(n_{j-1}),\cdots,\omega (N)),$
and $A_1>\cdots>A_k>0.$

For any decreasing sequence of $k$ natural numbers $a_1>\cdots>a_k>0,$ denote by $\mathcal{H}_N[a_1,\cdots,a_k]$ the set of $N$-step half-space walks $\omega$ such that
$\omega(0)=a,\ A_1(\omega)=a_1,\cdots, A_k(\omega)=a_k,\ n_k(\omega)=N.$
Particularly, $\mathcal{H}_N(\ell)$ is the set of $N$-step bridges of span $\ell$ for any $\ell\geq 0.$
Given an $\omega\in \mathcal{H}_N[a_1,\cdots,a_k],$ define the following new $N$-step walk $\omega^{\prime}:$ When $0\leq s\leq n_1(\omega),$
$\omega^{\prime}(s)=\omega(s).$ And when $s=n_1(\omega)+1,$
\begin{eqnarray*}
\omega^{\prime}(s)=\omega(s-1)\omega(s)^{-1}\omega(s-1)=\omega^{\prime}(s-1)\omega(s)^{-1}\omega(s-1).
\end{eqnarray*}
And recursively, when $n_1(\omega)+1<s\leq N,$
$\omega^{\prime}(s)=\omega^{\prime}(s-1)\omega(s)^{-1}\omega(s-1).$

Since $\left(\omega^{\prime}(s)\right)_{n_1(\omega)\leq s\leq N}$ is a reflection of $(\omega(s))_{n_1(\omega)\leq s\leq N},$
we see $\left(\omega^{\prime}(s)\right)_{n_1(\omega)\leq s\leq N}$ is an SAW. While $\left(\omega^{\prime}(s)\right)_{0\leq s\leq n_1(\omega)}$
is also an SAW, to prove $\omega^{\prime}$ is an SAW when $k\geq 2,$ it suffices to check that there is no cycle containing $\omega(n_1(\omega))$
in $\omega^{\prime}.$ {\it Actually}, when $k\geq 2,$ by the definition of $n_j(\omega)$'s, $\omega(n_j(\omega))^{-1}\omega(n_j(\omega)+1)$ must be $s_1^{-1},$
which implies that
$$\omega^{\prime}(n_1(\omega)+1)=\omega(n_1(\omega))s_1.$$
By our assumption, $s_1$ is an infinite-order generator independent of other generators $s_j\ (2\leq j\leq\ell_m),$ so there is no cycle
containing $\omega(n_1(\omega))$ in $\omega^{\prime}.$ Hence we have that
$$\omega^{\prime}\ \mbox{is an SAW and further}\ \omega^{\prime}\in \mathcal{H}_N[a_1+a_2,a_3,\cdots,a_k].\eqno{(2.2)}$$

Note that when crossing an edge of an SAW, increment of height function along this SAW is in $\{1,-1,0\};$ and Step 1.
It is easy to see that
$$\omega\in \mathcal{H}_N[a_1,\cdots,a_k]\rightarrow \omega^{\prime}\in \mathcal{H}_N[a_1+a_2,a_3,\cdots,a_k]\ \mbox{is an injective map.}$$
Thus $\vert \mathcal{H}_N[a_1,\cdots,a_k]\vert\leq\vert \mathcal{H}_N[a_1+a_2,a_3,\cdots,a_k]\vert\leq\cdots \leq\vert \mathcal{H}_N[a_1+\cdots+a_k]\vert,$ and further
\begin{eqnarray*}
&&h_N=\sum\left\vert \mathcal{H}_N[a_1,\cdots,a_k]\right\vert \\
&&\ \ \ \  \ \leq \sum\left\vert \mathcal{H}_N[a_1+\cdots +a_k]\right\vert =\sum b_{N,a_1+\cdots+a_k}=\sum\limits_{A=1}^NP_D(A)b_{N,A}
        \leq  P_D(N)b_N.
\end{eqnarray*}

\noindent{\bf Step 3.} Similarly to Proposition 6 in \cite{LiZY1}, for any constant $B>\pi \sqrt{\frac{2}{3}},$ there is an $N_0(B)$ satisfying
$$c_N\leq e^{BN^{\frac{1}{2}}}b_{N+1},\ \forall N\geq N_0(B).\eqno(2.3)$$
Note (2.3) holds for any $H_m$ and $H.$ Then for these graphs, the connective constant for bridges is just that for SAWs.
Let $b_n^{(m)}$ (resp. $b_n$) be the number of $n$-step bridges in $H_m$ (resp. $H$) starting at $v_m$
(resp. $v$), where $v_m$ is the induced vertex in $H_m$ of $v.$ Given any $\epsilon\in (0,1),$ for large enough $m,$
$$b^{(m)}_{\widehat{g}_m-1}=b_{\widehat{g}_m-1}\geq (\mu-\epsilon)^{\widehat{g}_m-1};$$
and further for any $s\geq 1,$
$b^{(m)}_{(\widehat{g}_m-1)s}\geq \left\{b^{(m)}_{\widehat{g}_m-1}\right\}^s=b^s_{\widehat{g}_m-1}\geq (\mu-\epsilon)^{(\widehat{g}_m-1)s},$
which implies that $\liminf\limits_{m\rightarrow\infty}\mu(H_m)\geq \mu(H).$ Clearly, $\mu(H_m)\leq \mu(H),\ m\geq 1.$
Hence $\lim\limits_{m\rightarrow\infty}\mu(H_m)=\mu(H).$\qed
\\

Recall marked groups were introduced in \cite{Gri}, and used to prove locality of percolation for abelian Cayley graphs in \cite{Mar-Tas}.
Here we extend a property of marked abelian groups in \cite{Mar-Tas} to marked finitely generated groups.

Let $d$ be a natural number. A $d$-marked finitely generated group is the data of finitely generated group $H$ with a
generating set $(s_1,s_2,\cdots, s_d)$, up to isomorphisms. And denote it as $[H; s_1,s_2,\cdots, s_d]$ or $H^{\bullet}$, depending on wether we want to point out the generating set or not. Here $[H_1; s_1,s_2,\cdots, s_d]$ and $[H_2;s'_1,s'_2,\cdots, s'_d]$ are isomorphic if there exists a group isomorphism from $H_1$ to $H_2$ mapping $s_i$ to $s_i'$ for all $i$.
Let $\mathbf{G}_d$ be the set of $d$-marked finitely generated groups.

Given a marked finitely generated group $H^{\bullet}=[H; s_1,s_2,\cdots, s_d]$ and a normal subgroup $\Lambda$ of $H$, the quotient $H^{\bullet}/\Lambda$ is denoted by
\[H^{\bullet}/\Lambda=\left[H/\Lambda; \overline{s_1},\overline{s_2}, \cdots,\overline{s_d}\right],\]
where $\left(\overline{s_1},\overline{s_2}, \cdots,\overline{s_d}\right)$ is the canonical image of $(s_1,s_2,\cdots, s_d)$.

Let $\delta=(\delta_1,\delta_2,\cdots,\delta_d)$ be the generating set of free group $F_d$. Recall that a finitely generated group $H$ with $d$ generators is isomorphic to a quotient group of free group $F_d$ by a normal subgroup $K$. Therefore, for any $d$-marked finitely generated group $H^{\bullet}=[H; s_1,s_2,\cdots, s_d]$, there is a unique normal subgroup $K$ of $F_d$ such that
\[H^{\bullet}\cong[F_d; \delta]/K=\left[F_d/K; \overline{\delta_1},\overline{\delta_2},\cdots, \overline{\delta_d}\right]. \eqno(2.4)\]
The  uniqueness of $K$ can be proved as follows. If there is another normal subgroup $K'$ such that
\[H^{\bullet}\cong [F_d; \delta]/K'=\left[F_d/K'; \overline{\delta_1}',\overline{\delta_2}',\cdots, \overline{\delta_d}'\right],\]
then there exists an isomorphism $\varphi$ from $F_d/K$ to $F_d/K'$ satisfying that $\varphi\left(\overline{\delta_i}\right)=\overline{\delta_i}'$, $1\leq i\leq d$.
Therefore, for any $t\in F_d$, $\varphi\left(\overline{t}\right)=\overline{t}'$, which forces $K=K'$.

By (2.4), $H^{\bullet}$ can be viewed naturally as a subset of $F_d,$ i.e., an element of $\{0,1\}^{F_d};$ and hence $\mathbf{G}_d$ can be viewed as a subset of $\{0,1\}^{F_d}.$ Endow $\mathbf{G}_d$ with the topology induced by the product topology on $\{0,1\}^{F_d}$. Then $\mathbf{G}_d$ is a Hausdorff compact space.
Let $\mathbf{G}$ be the set of all marked finitely generated groups, namely $\mathbf{G}$ is disjoint union of all the $\mathbf{G}_d$'s. Equip $\mathbf{G}$ with the topology generated by all open subsets of all the $\mathbf{G}_d$'s. {\bf Here and hereafter}, for any group $H$, $1_H$ is its identity element.

\begin{lem}\label{lem2}
Let $\{H_n^{\bullet}\}_{n=1}^\infty\subseteq \mathbf{G}$ be a sequence of marked finitely generated groups which converges to a finitely presented group $H^{\bullet}\in \mathbf{G}.$ Then $H^{\bullet}_n\cong H^{\bullet}/\Lambda_n$ for some subgroup $\Lambda_n$ of $H$ when $n$ is large enough; and for any fixed natural number $\ell$,
$\Lambda_n\cap B_H(1_H, \ell)=\{1_H\}$ for sufficiently large $n;$ and the relative girth of $H^{\bullet}/\Lambda_n$ to $H^{\bullet}$ tends to infinity. Particularly,
the corresponding Cayley graph sequence for $\{H_n^{\bullet}\}_{n=1}^\infty$ converges locally to the Cayley graph of $H^{\bullet}.$
\end{lem}
\pf
Assume $H^{\bullet}\in \mathbf{G}_d$ for some natural number $d$. By the assumption of the lemma, for large enough $n$, $H^{\bullet}_n\in \mathbf{G}_d$. By (2.4), for $n$ large enough, we have that
\[H^{\bullet}_n\cong [F_d; \delta]/K_n\ \mbox{and}\ H^{\bullet}\cong [F_d; \delta]/K.\]
Since $H^{\bullet}$ is finitely presented, we see that $K$ is finitely generated. Note $\{H_n^{\bullet}\}_{n=1}^\infty$ converges to $H^{\bullet}$ in $\mathbf{G}.$
Then for large enough $n$, a finite generating set $S$ of $K$ must be contained in $K_n$ and further $K$ is a subgroup of $K_n.$ Hence when $n$ is sufficiently large, let $\Lambda_n=K_n/K,$ we obtain $H^{\bullet}_n \cong H^{\bullet}/\Lambda_n.$

Since $H^{\bullet}/\Lambda_n$ converges to $H^{\bullet}$ in $\mathbf{G}_d,$ for any fixed natural number $\ell$, we have that
\[K_n\cap B_{F_d}(1_{F_d}, \ell)=K\cap B_{F_d}(1_{F_d}, \ell)\ \mbox{for large enough}\ n.\]
Thus when $n$ is sufficiently large,
\begin{align*}
  B_{F_d/K}(1_{F_d/K},\ell)\cap \Lambda_n=\Psi \left( B_{F_d}(1_{F_d},\ell)\cap K_n\right)=\Psi \left( B_{F_d}(1_{F_d},\ell)\cap K\right)=\{1_{F_d/K}\},
\end{align*}
where $\Psi: F_d\rightarrow F_d/K$ is the canonical quotient map.
Clearly, this implies that the relative girth of $H^{\bullet}/\Lambda_n$ to $H^{\bullet}$ tends to infinity as $n\rightarrow\infty.$
Therefore, the corresponding Cayley graph sequence of $\{H_n^{\bullet}\}_{n=1}^\infty$ converges locally to the Cayley graph of $H^{\bullet}.$
\qed

\begin{lem}\label{lem3}
For $\{G_n\}_{n=1}^\infty$ and $G$ specified in Theorem \ref{Thm1},
\[\lim_{n\rightarrow \infty}\mu(G_n)= \mu(G).\]
\end{lem}
\pf
Assume that  $\mu(G_n)\not \rightarrow \mu(G)$ as $n\rightarrow \infty$.
Then $\exists \epsilon>0$ such that
\[\limsup\limits_{n\rightarrow\infty} |\mu(G_n)- \mu(G)|>\epsilon. \eqno(2.5)\]

Let $G^{\bullet}_n$ (resp. $G^{\bullet}$) be the corresponding marked group for $G_n$ (resp. $G$). When (2.5) holds, without loss of generality, suppose $|\mu(G_n)- \mu(G)|>\epsilon$ for any $n\geq 1$ (otherwise choose a suitable subsequence). Since $\{G_n\}_{n=1}^\infty$ converges locally to $G$, we have that
for some natural number $d,$ $G^{\bullet}\in\mathbf{G}_d$ and $G^{\bullet}_n\in \mathbf{G}_d$ when $n$ is large enough.
For simplicity, we assume $G^{\bullet}_n\in \mathbf{G}_d$ for any $n\geq 1.$
Due to $\mathbf{G}_d$ is compact, we see that for some subsequence $\{n_k\}_{k=1}^\infty$ of natural numbers and $\widehat{G}^{\bullet}\in\mathbf{G}_d,$
\[G^{\bullet}_{n_k} \rightarrow \widehat{G}^{\bullet} ~\mbox{in} ~\mathbf{G}_d ~\mbox{as} ~k\rightarrow \infty.\]
By Lemma \ref{lem2}, we obtain that
$\widehat{G}^{\bullet}\cong G^{\bullet} ~\mbox{and}~ G^{\bullet}_{n_k} \cong  G^{\bullet}/\Lambda_{n_k}~\mbox{for sufficiently large}~k,$
where $\Lambda_{n_k}$ is some normal subgroup of $G^{\bullet}.$ And for any natural number $\ell,$ for large enough $k,$
$$\Lambda_{n_k}\cap B_{G^{\bullet}}(1_{G^{\bullet}}, \ell)=\{1_{G^{\bullet}}\},$$
and the relative girth of $G^{\bullet}/\Lambda_{n_k}$ to $G^{\bullet}$ tends to infinity as $k\rightarrow\infty.$

Now by Lemma \ref{lem1},
$\lim\limits_{k\rightarrow\infty}\mu(G_{n_k})= \mu(G).$
This is a contradiction to (2.5).\qed\\

{\bf  So far we have completed proving Theorem \ref{Thm1}.}

\begin{remark}\label{rem3}
From our proof, the following result holds: Assume $G\in\mathbf{\cal G}$ is a Cayley graph of a finitely generated group $\Gamma,$ and each $G_n\in\mathbf{\cal G}$ is a Cayley graph of a quotient group of $\Gamma.$ Then when $G_n$ converges locally to $G,$ $\lim\limits_{n\rightarrow\infty}\mu(G_n)=\mu(G).$
\end{remark}

\end{document}